\documentstyle[11pt]{article}

\begin{document}
\noindent
\begin{center}
  {\LARGE Discrete torsion and twisted orbifold cohomology}
  \end{center}

  \noindent
  \begin{center}

    {\large  Yongbin Ruan}\footnote{ partially
supported by the National Science Foundation }\\[5pt]
      Department of Mathematics, University of Wisconsin-Madison\\
        Madison, WI 53706\\[5pt]

              \end{center}

              \def \x{{\bf x}}
              \def \M{{\cal M}}
              \def \C{{\bf C}}
              \def \Z{{\bf Z}}
              \def \R{{\bf R}}
              \def \Q{{\bf Q}}
              \def \U{{\cal U}}
              \def \E{{\cal E}}
              \def \z{{\bf z}}
              \def \m{{\bf m}}
              \def \n{{\bf n}}
              \def \g{{\bf g}}
              \def \h{{\bf h}}
              \def \V{{\cal V}}
              \def \W{{\cal W}}
              \def \T{{\cal T}}
              \def \P{{\bf P}}
        \def \L{{\cal L}}

\section{Introduction}
    One of the remarkable insights of orbifold string theory is
    an indication of the existence of a new cohomology theory of
    orbifolds containing so-called twisted sectors as the contribution of singularities.
    Mathematically, such an orbifold cohomology theory has been
    constructed by Chen-Ruan \cite{CR}. Author believes that there
    is a "stringy" geometry and topology of orbifolds of which
    orbifold cohomology is its core. One aspect of this new
    geometry and topology is the twisted orbifold cohomology and its relation to discrete torison.
    Let me first explain their physical origin. Physicists usually work over a global quotient  $X=Y/G$ only, where $G$ is a finite group acting smoothly on
    $Y$. A discrete torsion is a cohomology class
    $\alpha\in H^2(G, U(1))$. Physically, a discrete torsion counts the freedom to choose a
    phase factor to weight path integral over each twisted sector without destroying the
    consistency of string
    theory.   For each $\alpha$, Vafa-Witten \cite{VW}  constructed the
    twisted orbifold cohomology group $H^*_{orb, \alpha}(X/G,
    \C)$.

    Mathematically, Vafa-Witten suggested that discrete torsion and twisted orbifold
    cohomology is connected to the problem of desingularizations. Recall that there are two
    methods to remove singularities, resolution or deformation. Both play important roles
    in the theory of Calabi-Yau 3-folds. One can
    obtain a smooth manifold $Y$ from an orbifold $X$ by using a combination of  resolution and
    deformation. We call $Y$ a desingularization of $X$. In string theory, we also require
    the resolution to be a crepant resolution.  It is known that a desingualization may not
    exist in dimension higher than three. In this case, we allow our desingularization to be
    an orbifold. As we mentioned in \cite{CR}, physicists
    predicted that ordinary orbifold cohomology group is the same as ordinary
    cohomology group  of its crepant resolution.  Vafa-Witten suggested that discrete
    torsion is a parameter for deformation. Furthermore, the cohomology of the
    desingularization is the twisted orbifold cohomology of discrete torsion plus possible
    contributions of exceptional loci of small resolution. A small resolution is a special
    kind of resolution whose exceptional loci is of codimension 2 or more. However, this proposal
    immediately ran into trouble because there are many more desingularizations than the
    number of discrete
    torsions. For example, D. Joyce \cite{JO} constructed five different desingularizations of
    $T^6/\Z_4$ while $H^2(\Z_4, U(1))=0$. To count these "missing" desingularizations seems to be
    a serious problem. On the another hand, it is well-known that  most  orbifolds (even Calabi-Yau orbifolds)
    are not global quotients.
    Therefore, it is also necessary to develop the theory over general orbifolds.

    We will address both problems in this paper. First, we introduce the notion
    of inner local system $\L$ for arbitrary orbifold. A local system is defined as an assignment of a flat
    (orbifold)-line bundle $L_{(g)}$ to each sector $X_{(g)}$ satisfying certain compatibility condition (See definition 2.1). Such a compatibility condition is
designed in such a way that Poincare duality and cup product in ordinary
orbifold cohomology survive the process of twisting.  Then, twisted orbifold
    cohomology $H^*_{orb}(X, \L)$ is defined as orbifold cohomology with value in the inner local system (See Definition 2.2).
    We will demonstrate that our inner local systems count all the examples
    constructed by D. Joyce. The author believes that the inner local system is a more
    fundamental notion
    than the discrete torsion. Then, we can formulate
    following mathematical conjecture: {\em Suppose that $X$ is a Calabi-Yau
    Gorenstein
    orbifold. For every desingularization, we can associate a
    inner local systems such that as additive groups the ordinary   orbifold cohomology
    of desigularization is the  sum of  twisted orbifold cohomology and contributions from exceptional loci of small resolution.}

     Our next goal is to determine appropriate notion of discrete torsion for general
    orbifold.  Let $X$ be an
    arbitrary almost complex orbifold. The author's key observation is that we should use the orbifold
    fundamental group $\pi^{orb}_1(X)$ (See definition 2.1) to replace $G$. Then a
    discrete torsion of $X$ is defined as a cohomology class
    $\alpha\in H^2(\pi^{orb}_1(X), U(1))$. Note that if
    $X=Y/G$ is a global quotient, there is a short exact sequence
    $$1\rightarrow \pi_1(Y)\rightarrow \pi^{orb}_1(X)\rightarrow G\rightarrow 1. \leqno(1.1)$$
     It induces a
    homomorphism $H^2(G, U(1))\rightarrow H^2(\pi^{orb}_1(X),
    U(1))$. Hence a discrete torsion in the sense of Vafa-Witten
    induces the discrete torsion in this paper. In fact, we can do better, we can
    define a local discrete torsion for each connected component of singular loci.
    A global discrete torsion is defined as an assignment of a local discrete torsion to
    a connected component of singular loci. Then, the link between discrete
    torsion and twisted orbifold cohomology is the theorem
    that a global discrete torsion induces an inner local system and hence
    define a twisted orbifold cohomology. However, we want to emphasis that not every inner local
    system comes from discrete torsion (See example 5.3).

    We will introduce inner local system and twisted orbifold cohomology ring in
    section 2. The section 3 is devoted to discrete torsion.
    The relation between discrete torsion and local system is discussed in
    section 4. Finally, some examples are computed
    in section 5. This paper can be viewed as a sequel to \cite{CR}.
    Since many constructions are similar, we will follow the notations
    of \cite{CR} and be sketchy in the details.
    The author strongly encourages readers to read \cite{CR} first before
    reading this paper.

    This paper was completed while author was visiting Caltech. He would like
    to thank R. Pandharipande and the Caltech math department for financial  support and
    hospitality. The author would like to thank E. Zaslow
    for bringing  his
    attention to Vafa-Witten's paper and E. Witten for explaining to
    him \cite{VW}. He would like also thank A. Adem for many
    valuable discussions about group cohomology.

\section{Local system and twisted orbifold cohomology}

\subsection{Review of ordinary orbifold cohomology}

        Suppose that $X$ is an orbifold. By the definition, $X$ is a topological space
    with a system of orbifold charts (uniformizing system). Namely,  every point
        $p\in X$ has a system of orbifold chart  of the form
        $U_p/ G_p$ where $U_p$ is a smooth manifold  and
        $G_p$ is a finite group acting on $U_p$ fixing $p$.  $G_p$ is called a local group.
        Note that the action of $G_p$ does not have to be effective. If it does,
        we call it a reduced orbifold. We use $(U_p, G_p)$
    to denote the chart.
The patching condition is
        follows: if $q\in U_p/G_p\cap U_r/G_r$, there is an orbifold chart $U_q/G_q\subset
    U_p/G_p\cap U_r/G_r$. Moreover, the inclusion map
        $i: U_q/G_q\subset U_p/G_p$  can be lifted to a smooth map
        $$\tilde{i}_{pq}: U_q\rightarrow U_p \leqno(2.1)$$
        and an injective homomorphism
        $$i_{\#,pq}: G_q\rightarrow G_p\leqno(2.2)$$
        such that $\tilde{i}$ is $i_{\#}$-equivariant.
        $$i_{pq}=(\tilde{i},i_{\#}): (U_q,G_q)\rightarrow (U_p, G_p)$$
    is called an injection.
        A different lifting differs from $\tilde{i}$ by the action of an element of $G_p$.
        Moreover, $i_{\#}$ differs by the conjugation of the same element. We say that the corresponding
    injections are equivalent.

    Therefore, for any $g\in G_q$, the conjugacy class $(i_{\#}(g))_{G_p}$ is
    well-defined. We define an equivalence relation $(g)_{G_q}\cong
    (i_{\#}(g))_{G_p}$. Let $T_1$ be the set of equivalence classes. By abusing the
    notation, we use $(g)$ to denote the equivalence class to which
    $(g)_{G_q}$     belongs to. For each $(g)\in T_1$, we can define a sector
    $$X_{(g)}=\{(x, (g')_{G_x})| g'\in G_x, (g')_{G_x}\in (g)\}.\leqno(2.3)$$
    It was shown in \cite{CR} that $X_{(g)}$ is an orbifold. It is the common
    convention that we call $X_{(g)}$ for $g\neq 1$ a twisted sector
    and $X_{(1)}$ a nontwisted sector. Once we define sectors, we
     diagonalize the action of $g$ in  $T_x X_{(g)}$
     for each $x\in X_{(g)}$. Suppose that
    $g=diag(e^{\frac{2\pi i m_1}{m}}, \cdots, e^{\frac{2\pi im_n}{m}})$, where
    $m$ is the order of $g$ and $0\leq \frac{m_i}{m}<1$.
    Then we define the degree shifting number $\iota_{(g)}=\sum_i
    \frac{m_i}{m}$. One can show that $\iota_{(g)}$ is independent of $x\in
    X_{(g)}$. The ordinary orbifold cohomology is defined as
    $$H^*_{orb}(X, \C)=\oplus_{(g)\in T_1} H^{*-2\iota_{(g)}}(X_{(g)}, \C).
    \leqno(2.4)$$
    There is a diffeomorphism $I: X_{(g)}\rightarrow X_{(g^{-1)}}$
    defined by $(x, (g))\rightarrow (x, (g^{-1}))$.
    Poincare paring $<\,\ >_{orb}$ of orbifold cohomology is defined as the direct
    sum of
    $$<\, \ >^{(g)}_{orb}: H^{d- 2\iota_{(g)}}
    (X_{(g)},\C)\otimes
    H^{2n-d-2\iota_{(g^{-1})}}(X_{(g^{-1})}, \C)\rightarrow \C\leqno(2.5)$$
    defined by
    $$<\alpha, \beta>^{(g)}_{orb}=\int_{X_{(g)}} \alpha\wedge I^*\beta.\leqno(2.6)$$

    Next, we consider cup product.  We first construct a moduli space
    (see \cite{CR} (section 4.1)).
    $$X_3=\{(x, (g_1,g_2,g_3)_{G_x})| g_i\in G_x, g_1g_2g_3=1\}\leqno(2.7)$$
    $X_3$ is an orbifold.
    Let $\g=(g_1,\cdots, g_k)$ with  $g_i\in G_q$. By abusing
    the notation, we simply say $\g\in G_q$.  We define the equivalence relation $(\g)_{G_q}\cong
    (i_{\#}(\g))_{G_q}$. Let $T_k$ be the set  of equivalence class and use
    $(\g)$ to denote the equivalence class such that $(\g)_{G_q}\in (\g)$.
    We will use $T^o_k\subset T_k$ to denote the set of equivalence classes
    of $(\g)$ such that  $\g=(g_1, \cdots, g_k)$ with $g_1\cdots g_k=1$.
    It was proved in \cite{CR} that
    $$X_{(\g)}=\{(x, (\g')_{G_x})| \g'\in G_x, (\g')=(\g)\}.\leqno(2.8)$$
    is an orbifold.
    One can check that
    $$X_3=\bigsqcup_{(\g)\in T^o_3} X_{(\g)},\leqno(2.9)$$
Then, for each $(\g)\in T^o_3$, we can define evaluation maps
    $$e_i: X_{(\g)} \rightarrow X_{(g)}$$
    by
    $$e_i(x, (g_1,g_2,g_3)_{G_x})=(x, (g_i)_{G_x}).$$
    Furthermore, there is an obstruction bundle $E$ (see  Lemma 4.2.2 \cite{CR}).
    Then, we can define a three-point function
    $$<\alpha, \beta, \gamma>_{orb}=\int_{X_{(\g)}} e^*_1\alpha \wedge
    e^*_2\beta\wedge e^*_3\gamma\wedge e(E),\leqno(2.10)$$
    for any $\alpha\in H^{p-\iota_{(g_1)}}(X_{(g_1)}, \C),  \beta\in H^{q-\iota_{(g_2)}}
    (X_{(g_2)}, \C),\gamma\in H^{2n-p-q-\iota_{(g_3)}}(X_{(g_3)}, \C)$.
    Once the three point function is defined, the cup product is defined
    by the equation
    $$<\alpha\cup_{orb}\beta, \gamma>_{orb}=<\alpha, \beta, \gamma>_{orb}.
    \leqno(2.11)$$
    for arbitrary $\gamma$.

    There is a Dolbeaut version of orbifold cohomology ring (Dolbeaut orbifold
    cohomology ring) with identical construction. We refer reader to
    \cite{CR} for detail.

    Next, we consider the bundle over each sector and its pull-back. By \cite{CR1},
    this is a very subtle problem and one has to be careful. We will give explicit
    construction in our case and refer reader to general theory in \cite{CR1}.

    Now, let's examine the orbifold structure of twisted sectors more carefully.
    Suppose that $p\in X_{(g)}$ and an orbifold chart of $p\in X$ is $U_p/G_p$.
    By \cite{CR1}(Lemma 3.1.1), an orbifold chart of $X_{(g)}$ can be described as
    follows. Choose a representative of $(g)_{G_p}$, say $g_p$. Then, a local orbifold
    chart of $p$ is $U_{g_p}/C(g_p)$, where $U_{g_p}$ is the fixed point loci of $g_p$ and
    $C(g_p)$ is the centralizer. In general, $C(g_p)$ may not acts freely on $U_{g_p}$.
     The patching map of $X_{(g)}$ is defined
    in the same way. More generally,    $X_{(\g)}$ also has a structure of a orbifold given by
    $(U_{\g_p},C(\g_p))$, where $U_{\g_p}, C(\g_p)$ are the fixed point loci and centralizer of
      $\g_p$. We denote the corresponding patching map by injection
    $i_{\g,pq}=(\tilde{i}_{pq, (\g)}, i_{\#, pq, (\g)})$.

    Now, an orbifold-bundle $f: E\rightarrow X_{(\g)}$ is a continuous map between topological space
    such that $E$ has a structure of orbifold as follows. Suppose that $p\in X_{(\g)}$. $E$ is covered
    by chart of the form $(U_{\g_p}\times \R^n, C(\g_p))$ such that the restriction of $f$ is
    the projection  $U_{\g_p}\times \R^n\rightarrow U_{\g_p}$ equivariant under $C(\g_p)$. For any injection
    $i: (U_{\g_q}, C(\g_q))\rightarrow (U_{\g_p},C(\g_p))$, there is an injection
    between charts $(U_{\g_q}\times \R^n, G_{\g_q})\rightarrow (U_{\g_p}\times \R^n, G_{\g_p})$
    given by $(\tilde{i}_{pq, (\g)}\times g_i, i_{\#,pq, (\g)}),$ where $g_{i}: U_{\g_q}\rightarrow
    Aut(\R^n)$ satisfies the condition $g_{j\circ i}=g_j\circ g_i$. The last condition is to ensure
     that the equivalent injections on $X_{(\g)}$  implies the equivalence between corresponding
    injections on total space of bundle.

      Next, we consider evaluation map $e_i: X_{(\g)}\rightarrow X_{(g)}$. A crucial observation is
    that we can choose the orbifold charts of $X_{(\g)}, X_{(g)}$ silmontaneously. Suppose that
    $p\in X_{(\g)}$. We first  choose $\g_p$. Then, it gives
    a natural choice for $g_{i,p}$. Similarly, an injection $i_{pq, (\g)}$ between the charts of
    $X_{(\g)}$  gives a natural choice of injection $\lambda(i_{pq, (\g)})=i_{pq, (g_i)}$ with the
    property $\lambda(j\circ i)=\lambda(j)\circ \lambda(i)$.
    The evaluation map is interpreted as an inclusion
    $$e_{i,p}: U_{\g_p}\rightarrow U_{g_{i,p}}, e_{\#,i,p}: C(\g_p)
    \rightarrow C(g_{i,p}).\leqno(2.11)$$
    Following \cite{CR1}, we say that $e_i$ is a good map. By \cite{CR1}, if $E\rightarrow X_{(g_i)}$
    is a orbifold-bundle, $e^*E\rightarrow X_{(\g)}$ is a orbifold-bundle.

    Note that the direct sum and tensor product of orbifold-bundles is still a orbifold bundle. Moreover, all
    the differential geometric constructions such as differential form, connection and curvature work over
    orbifold-bundle.

\subsection{Inner Local system and twisted orbifold cohomology ring}
    Now, we introduce the notion of inner local system for orbifold.
    \vskip 0.1in
    \noindent
    {\bf Definition 2.1: }{\it Suppose that $X$ is an orbifold (almost complex or
    not). An inner  local system $\L=\{L_{(g)}\}_{g\in T_1}$ is an assignment of  a flat
    complex orbifold line bundle over
    $$L_{(g)}\rightarrow X_{(g)}$$
    to each sector
    $X_{(g)}$ satisfying the compatibility condition
    \begin{description}
    \item[(1)] $L_{(1)}=1$ is trivial.
    \item[(2)] $I^*L_{(g^{-1})}=L^{-1}_{(g)}.$
    \item[(3)] Over each $X_{(\g)}$ with $(\g)\in T^o_3$,
    $\otimes_i e^*_i L_{(g_i)}=1 $.
    \end{description}
    If $X$ is a complex orbifold, we assume that $L_{(g)}$ is holomorphic.}
    \vskip 0.1in
        \vskip 0.1in
    \noindent
    {\bf Definition 2.2: }{\it Given an inner local system $\L$, we define the twisted orbifold cohomology
    $$H^*_{orb}(X, \L)=\oplus_{(g)}
    H^{*-2\iota_{(g)}}(X_{(g)}, L_{(g)}).$$}
    \vskip 0.1in
    \noindent
    {\bf Definition 2.3: }{\it Suppose that $X$ is a closed complex
    orbifold and $\L$ is an inner local system. We define Dolbeault cohomology groups
    $$H^{p,q}_{orb}(X, \L)=\oplus_{(g)}H^{p-\iota_{(g)},
    q-\iota_{(g)}}(X_{(g)}; L_{(g)}).\leqno(2.12)$$}
    \vskip 0.1in
    \noindent
    {\bf Proposition  2.4: }{\it If $X$ is a Kahler orbifold, we have
    Hodge decomposition
    $$H^k_{orb}(X, \L)=
\oplus_{k=p+q} H^{p,q}_{orb}(X, \L).\leqno(2.13)$$}
    \vskip 0.1in
    \noindent
    {\bf Proof: } Note that each sector $X_{(g)}$ is a K\"ahler orbifold.
    The proposition follows by applying the ordinary Hodge theorem with twisted
    coefficients to each sector $X_{(g)}$. $\Box$

    The property (2) of Definition 2.1 can be used to show that Poincare pairing defined in (2.6)
    can be adopted to twisted orbifold cohomology.
 \vskip 0.1in
    \noindent
    {\bf Definition (Poincar\'e duality) 2.5: }{\it Suppose that $X$ is a $2n$-dimensional
    closed almost complex orbifold. We define a pairing
    $$<\, \ >_{orb, \L}: H^{d}_{orb}(X,\L)\otimes
    H^{2n-d}_{orb}(X, \L)\rightarrow \C.\leqno(2.14)$$
    as the direct sum of
    $$<\, \ >^{(g)}_{orb, \L}: H^{d- 2\iota_{(g)}}
    (X_{(g)},L_{(g)})\otimes
    H^{2n-d-2\iota_{(g^{-1})}}(X_{(g^{-1})}, L_{(g^{-1})})\rightarrow \C\leqno(2.15)$$
    defined by
    $$<\alpha, \beta>^{(g)}_{orb, \L}=\int_{X_{(g)}} \alpha\wedge I^*\beta.\leqno(2.16)$$
    }
    \vskip 0.1in
    Note that $L_{(g)}I^*L_{(g^{-1})}=1$. Hence the integral (2.6) makes
    sense.
    \vskip 0.1in
    \noindent
    {\bf Theorem 2.6: }{\it The pairing $<\,\ >_{orb,\L}$ is
    nondegenerate.}
    \vskip 0.1in
    \noindent
    {\bf Proof: } The proof follows from ordinary Poincare duality on
    $X_{(g)}$ with twisted coefficient.

    There is also a version of Poincar\'e duality for twisted Dolbeault
    cohomology. Suppose that $X$ is a
    closed  complex orbifold of complex dimension $n$. Then    $X_{(g)}$   is a closed complex orbifold.
    \vskip 0.1in
    \noindent
    {\bf Definition 2.7: }{\it We define a pairing
    $$<\, \ >_{orb, \L}: H^{p,q}_{orb}(X,\L)\otimes
    H^{n-p, n-q}_{orb}(X, \L)\rightarrow \C.\leqno(2.17)$$
    as the direct sum of
    $$<\, \ >^{(g)}_{orb, \L}: H^{p-\iota_{(g)},q-\iota_{(g)}}
    (X_{(g)},L_{(g)})\otimes
    H^{n-p-\iota_{(g^{-1})},n-q-\iota_{(g^{-1})}}(X_{(g^{-1})},
    L_{(g^{-1})})\rightarrow \C\leqno(2.18)$$
    defined by
    $$<\alpha, \beta>^{(g)}_{orb, \L}=\int_{X_{(g)}} \alpha\wedge I^*\beta.\leqno(2.19)$$
    }
    \vskip 0.1in
        \noindent
    {\bf Theorem 2.8: }{\it The  pairing (2.17) is
    nondegenerate.}
    \vskip 0.1in

    The property (3) of Definition 2.1 shows that the integral (2.10) makes sense for twisted
    orbifold cohomology classes. The same construction of \cite{CR} goes through
    without  change. We can define a twisted orbifold product
    $\cup_{orb, \L}$ in the same fashion. The same proof as in \cite{CR} yields
    \vspace{3mm}

\noindent{\bf Theorem 2.9: } {\it Let $X$ be a closed almost
complex orbifold with almost complex structure $J$ and $\dim_\C
X=n$. There is a cup product $\cup_{orb,\L}:
H_{orb}^p(X;\L)\times H_{orb}^q(X;\L) \rightarrow
H_{orb}^{p+q}(X;\L)$ for any $0\leq p,q\leq 2n$ such that
$p+q\leq 2n$, which has the following properties:
\begin{enumerate}
\item The total twisted orbifold cohomology group $H^\ast_{orb}(X;\L)=
\oplus_{0\leq d\leq 2n}H^d_{orb}(X;\L)$ is a ring with unit
$e_X^0\in H_{orb}^0(X;\L)$ under $\cup_{orb,\L}$, where $e_X^0$ is
the Poincare\'{e} dual to the fundamental class of nontwisted
sector.
\item Restricted to each $H_{orb}^d(X;\L)\times H_{orb}^{2n-d}(X;\L)
\rightarrow H_{orb}^{2n}(X;\L)=H^{2n}(X, \C)$, $$ \int_X
\alpha\cup_{orb,\L} \beta=<\alpha, \beta>_{orb,\L}.
\leqno(2.20) $$

\item The cup product $\cup_{orb,\L}$ is invariant under  deformations of
$J$.
\item When $X$ is of integral degree shifting numbers, the total twisted orbifold
cohomology group $H_{orb}^\ast(X;\L)$ is integrally graded, and we
have supercommutativity $$
\alpha_1\cup_{orb,\L}\alpha_2=(-1)^{\deg\alpha_1\cdot\deg\alpha_2}\alpha_2
\cup_{orb,\L}\alpha_1. $$
\item Restricted to the nontwisted sectors, i.e., the ordinary
cohomology $H^\ast(X;\C)$, the cup product $\cup_{orb,\L}$ equals the
ordinary cup product on $X$.
\item $\cup_{orb, \L}$ is associative.
\end{enumerate}}

\vspace{3mm}
    Similarly, we also have a holomorphic version.
    \vspace{3mm}

\noindent{\bf Theorem 2.10: }{\it Let $X$ be an n-dimensional closed
complex orbifold with complex structure $J$.  The orbifold cup
product $$ \cup_{orb,\L}: H_{orb}^{p,q}(X;\L)\times
H_{orb}^{p^\prime,q^\prime}(X;\L) \rightarrow
H_{orb}^{p+p^\prime,q+q^\prime}(X;\L) $$ has the following
properties:
\begin{enumerate}
\item The total orbifold Dolbeault cohomology group is a ring with unit
$e_X^0\in H_{orb}^{0,0}(X;\L)$ under $\cup_{orb,\L}$, where $e_X^0$
is the class represented by the equal one constant function on
$X$.
\item Restricted to each $H_{orb}^{p,q}(X;\L)\times H_{orb}^{n-p,n-q}(X;\L)
\rightarrow H_{orb}^{n,n}(X;\L)$, the integral
$\int_{X}\alpha\cup_{orb,\L} \beta$ equals the Poincare pairing
$<\alpha,\beta>_{orb,\L}$.
\item The cup product $\cup_{orb,\L}$ is invariant under the
deformation of $J$.
\item When $X$ is of integral degree shifting numbers, the total twisted orbifold
Dolbeault cohomology group of $X$ is integrally graded, and we
have supercommutativity $$
\alpha_1\cup_{orb,\L}\alpha_2=(-1)^{\deg\alpha_1\cdot\deg\alpha_2}\alpha_2
\cup_{orb,\L}\alpha_1. $$
\item Restricted to the nontwisted sectors, i.e., the ordinary Dolbeault
cohomologies $H^{\ast,\ast}(X;\C)$, the cup product $\cup_{orb,\L}$
equals the ordinary wedge product on $X$.
\item The cup product is associative.
\item When $X$ is K\"ahler, the cup product $\cup_{orb,\L}$ coincides with
the twisted orbifold cup product over the twisted orbifold
cohomology groups $H^\ast_{orb}(X;\L)$ under the relation $$
H^r_{orb}(X;\L)\otimes\C=\oplus_{p+q=r}H^{p,q}_{orb}(X;\L).
$$
\end{enumerate}}
\vskip 0.1in
    \noindent
    {\bf Remark 2.11: }{\it If $X$ is open, we can define usual
    twisted  orbifold cohomology $H^*_{orb}(X, \L)$ and twisted orbifold cohomology
    with compact support $H^*_{orb, c}(X, \L)$ in the same
    fashion. The Poincare paring should be understood as the
    paring between $H^{d}_{orb}(X, \L)$ and $H^{2n-d}_{orb, c}(X,
    \L)$.}\vskip 0.1in

\section{Orbifold fundamental group and discrete torsion}

    First, we recall the definition of orbifold fundamental group.
                \vskip 0.1in
        \noindent
        {\bf Definition 3.1: }{\it A smooth map $f:Y\rightarrow X$
        is an orbifold cover iff (1) each $p\in Y$ has a neighborhood
      $U_p/G_p$ such that the restriction of $f$ to $U_p/G_p$ is    isomorphic
      to a map $U_p/G_p\rightarrow U_p/\Gamma$ such that $G_p\subset
    \Gamma$ is a subgroup. (2) Each $q\in X$ has a neighborhood
    $U_q/G_q$ for which each component of $f^{-1}(U_q/G_q)$ is
    isomorphic to $U_q/\Gamma'$ such that $\Gamma'\subset G_q$ is a
    subgroup.           An
        orbifold universal cover $f:Y\rightarrow X$ of $X$ has the
        property: (i) $Y$ is connected; (ii)if $f': Y'\rightarrow X$ is an orbifold
        cover, then there exists an orbifold cover $h:
        Y\rightarrow Y'$ such that $f=f'\circ h$. If $Y$ exists, we call $Y$ the
        orbifold universal cover of $X$ and the group of deck translations the orbifold
        fundamental group $\pi^{orb}_1(X)$ of $X$.}
        \vskip 0.1in
        By Thurston \cite{T}, an orbifold universal cover exists.
         It is clear from the definition that the orbifold universal
        cover is unique. Suppose that $f: Y\rightarrow X$ is an orbifold universal
        cover. Then
        $$f: Y-f^{-1}(\Sigma X)\rightarrow X-\Sigma X\leqno(3.1)$$
        is an honest cover with $G=\pi^{orb}_1(X)$ as covering
        group, where $\Sigma$ is the singular loci of $X$. Therefore, $X=Y/G$ and there is a surjective
        homomorphism
        $$p_f: \pi_1(X-\Sigma X)\rightarrow G.\leqno(3.2)$$
        In general, (3.1) is not a universal covering. Hence,
        $p_f$ is not an isomorphism.

        \vskip 0.1in
        \noindent
        {\bf Remark 3.2: }{\it Suppose that $X=Z/G$ for an orbifold $Z$
     and $Y$ is the orbifold universal cover of $Z$. By the
        definition, $Y$ is an orbifold universal cover of $X$. It
        is clear that there is a short exact sequence
        $$1\rightarrow \pi_1(Z)\rightarrow \pi^{orb}(X)\rightarrow
        G\rightarrow 1.\leqno(3.3)$$}

        \vskip 0.1in
        \noindent
        {\bf Example 3.3: } Consider the Kummer surface $T^4/\tau$
        where $\tau$ is the involution
        $$\tau(e^{it_1}, e^{it_2}, e^{it_3}, e^{it_4})=
        (e^{-it_1}, e^{-it_2}, e^{-it_3}, e^{-it_4}).\leqno(3.4)$$
        The universal cover is $\R^4$. The group $G$ of deck
        translations is generated by translations $\lambda_i$ by an integral point
        and the involution
    $$\tau: (t_1, t_2, t_3, t_4)\rightarrow (-t_1,
        -t_2, -t_3, -t_4).$$
     It is easy to check that
        $$G=\{\lambda_i (i=1,2,3,4), \tau| \tau^2=1,
        \tau\lambda_i=\lambda^{-1}_i\tau,\}\leqno(3.5)$$
        where $\lambda_i$ represents translation and $\tau$ represents
          involution.
        \vskip 0.1in
        \noindent
        {\bf Example 3.4: } Let $T^6=\R^6/\Gamma$ where $\Gamma$ is the lattice
        of integral points. Consider $\Z^2_2$ acting on $T^6$ lifted
        to an action on $\R^6$ as
        $$\sigma_1(t_1, t_2, t_3, t_4, t_5, t_6)=(-t_1, -t_2, -t_3, -t_4, t_5,
        t_6)$$
        $$\sigma_2(t_1, t_2, t_3, t_4, t_5, t_6)=(-t_1, -t_2, t_3, t_4, -t_5,
        -t_6)$$
        $$\sigma_3(t_1, t_2, t_3, t_4, t_5, t_6)=(t_1, t_2, -t_3, -t_4, -t_5,
        -t_6).$$
        This  example was considered by Vafa-Witten \cite{VW}.
        The orbifold fundamental group
        $$\pi^{orb}_1(T^6/\Z^2_2)=\{\tau_i (1\leq i\leq 6), \sigma_j
        (1\leq j\leq 3)|$$
        $$\sigma^2_i=1, \sigma_1\tau_i=\tau^{-1}_i\sigma_1 (i\neq
        5,6), \sigma_2\tau_i=\tau^{-1}_i\sigma_2 (i\neq 3,4),
        \sigma_3\tau_i=\tau^{-1}_i\sigma_3 (i\neq
        1,2)\}.\leqno(3.6)$$
        \vskip 0.1in
        The following example was taken from \cite{SC}
        \vskip 0.1in
        \noindent
        {\bf Example 3.5: }Consider the orbifold Riemann surface
        $\Sigma_g$ of genus $g$ and $n$-orbifold points $\z=(x_1, \cdots, x_n)$
        with orders $k_1, \cdots, k_n$. Then,
        $$\pi_1^{orb}(\Sigma_g)=\{\lambda_i (i\leq 2g), \sigma_i (i\leq n)|
        \sigma_1\cdots \sigma_n\prod_i [\lambda_{2i-1}, \lambda_{2i}]=1,
        \sigma^{k_i}_i=1\},\leqno(3.7)$$
        where $\lambda_i$ are the generators of $\pi_1(\Sigma_g)$ and
        $\sigma_i$ are the generators of $\Sigma_g-\z$ represented by a loop
        around each orbifold point.

         Note that $\pi_1^{orb}(\Sigma_g)$ is just  $\pi_1(\Sigma_g-\z)$ modulo
        by the relation $\sigma^{k_i}_i=1$. This suggests that one can first take the cover
        of $\Sigma_g-\z$ induced by $\pi^{orb}_1(\Sigma)$. The relation $\sigma^{k_i}_i=1$
        implies that the preimage of the punctured disc around $x_i$ is a punctured disc.
        Then we can fill in the center point to obtain the orbifold
        universal cover.
        \vskip 0.1in
        \noindent
        {\bf Definition 3.6: }{ Suppose that $S\subset X$ is a connected component of
    singular loci. A local discrete torsion $\alpha_S$ at $S$ is defined as a cohomology
    class $\alpha_S\in H^2(\pi^{orb}_1(U(S)), U(1))=H^2(B\pi^{orb}_1(U(S)), U(1))$, where
    $U(S)$ is a small open neighborhood of $S$. A global discrete torsion $\alpha=
    \{\alpha_S\}$ is an assignment of a local discrete torsion to each connected component
    of singular loci.  }
        \vskip 0.1in
        If $X=Z/G$ for a finite group $G$, by Remark 3.2, there is a surjective homomorphism
        $$\pi: \pi^{orb}_1(X)\rightarrow G.$$
        $\pi$ induces a homomorphism
        $$\pi^*: H^2(G, U(1))\rightarrow H^2(\pi^{orb}_1(X),
        U(1)).\leqno(3.8)$$
        Hence, an element of $H^2(G, U(1))$ induces a discrete
        torsion of $X$.

        They are many ways to define $H^2(G, U(1))$. The definition $H^2(G, U(1))=H^2(BG, U(1))$ is a very useful definition
    for computation since we can use algebro-topological machinery. However,
          we can also take the original definition in terms of cocycles.
    A 2-cocycle is a map $\alpha: G\times
         G\rightarrow U(1)$ satisfying
    $$\alpha_{g,1}=\alpha_{1,g}=1, \alpha_{g,hk}\alpha_{h,k}
         =\alpha_{g,h}\alpha_{gh,k}, \leqno(3.9)$$
         for any $g,h,k\in G$. We denote the set of two cocycles by
         $Z^2(G, U(1))$.
    For any map $\rho: G\rightarrow U(1)$ with $\rho_1=1$, its coboundary
    is defined by formula
    $$(\delta \rho)_{g,h}=\rho_g\rho_h\rho_{gh}^{-1}.\leqno(3.10)$$
    Let $B^2(G, U(1))$ be the set
    of coboundaries. Then, $H^2(G,U(1))=Z^2(G,U(1))/B^2(G,U(1))$.
    $H^2(G, U(1))$ naturally appears in many important places of mathematics.
    For example, it classifies the group extension of $G$ by $U(1)$. If
    we have a unitary projective representation of $G$, it naturally induces
    a class of $H^2(G, U(1))$. In many instances, this class completely classifies
    the projective unitary representation. In fact, it is in this context that
    discrete torsion arises in orbifold string theory.

         \noindent
         {\bf Definition 3.7: }{\it
         For each 2-cocycle $\alpha$, we define its phase
         $$\gamma(\alpha)_{g,h}=\alpha_{g,h}\alpha^{-1}_{h,g}.\leqno(3.11)$$
         }
         \vskip 0.1in
         It is clear that
         $\gamma(\alpha)_{g,g}=1,\gamma(\alpha)_{g,h}=\gamma(\alpha)^{-1}_{h,g}.$
         \vskip 0.1in
         \noindent
         {\bf Lemma 3.8: }{\it Suppose that $gh=hg, gk=kg.$ Then
         \begin{description}
          \item[(1)] $\gamma(\delta \rho)_{g,h}=1$.

          \item[(2)] $\gamma(\alpha)_{g,
          hk}=\gamma(\alpha)_{g,h}\gamma(\alpha)_{g,k}.$
          \end{description}
          The (2) implies $L^{\alpha}_g=\gamma_{g,\cdot}:
          C(g)\rightarrow U(1)$ is a group homomorphism. We call
          $L^{\alpha}_g$ a $\alpha$-twisted character.
          }
          \vskip 0.1in
          \noindent
          {\bf Proof: } (1) is obvious. For (2),
          $$\begin{array}{lll}
           \gamma(\alpha)_{g,hk}&=&\alpha_{g,hk}\alpha_{hk,g}^{-1}\\
                                &=&\alpha_{g,hk}\alpha^{-1}_{gh,k}\alpha_{hg,k}
                                \alpha_{h,gk}^{-1}\alpha_{h,kg}\alpha^{-1}_{hk,g}\\
                                &=&\alpha_{g,h}\alpha^{-1}_{h,k}\alpha_{g,k}
                                    \alpha^{-1}_{h,g}\alpha_{h,k}\alpha^{-1}_{k,g}\\
                                &=&\gamma(\alpha)_{g,h}\gamma(\alpha)_{g,k}
           \end{array}
           $$
           \vskip 0.1in
            Next, we calculate discrete torsions for some groups.
    We first consider the case of finite abelian group $G$. In this case
    $H^i(G, \Q)=0$ for $i\neq 0$. The exact sequence
    $$0\rightarrow \Z\rightarrow \C\rightarrow \C^*\rightarrow 1$$
    implies that $H^2(G,U(1))=H^2(G,\C^*)=H^3(G, \Z)$. By universal coefficient
    theorem, $H^3(G,\Z)= H_2(G,\Z).$
    \vskip 0.1in
           \noindent
           {\bf Example 3.9 $G=\Z/n\times \Z/m$: } Notes that $H^2(G, U(1))=H_2(G,\Z)=
    \Z/n\otimes \Z/m=Z_{gcd(n,m)}$. In this case, one can write down the phase of
    discrete torsion explicitly \cite{VW}. Let $\xi (\zeta)$ be $n(m)$-root
    of unity. Any element of $\Z/n\times \Z/m$ can be written as
    $(\xi^a, \zeta^b)$. Let $p=gcd(n,m)$. The phase of a discrete torsion can be
    written as
    $$\gamma_{(\xi^a, \zeta^b),(\xi^{a'},
    \zeta^{b'})}=\omega_p^{m(ab'-ba')}$$
    with $\omega_p=e^{2\pi i/p}, m=1, \cdots, p.$ There are
    $p$-phases for $p$-discrete torsions.  It is trivial
    to generalize this construction to an arbitrary finite abelian group.

\section{Discrete torsion and local system}

Suppose that $f: Y\rightarrow X$ is the orbifold universal
    cover and $G$ is the orbifold fundamental group which acts
    on $Y$ such that $X=Y/G$. Suppose $X_{(g)}$ is a
    sector (twisted or nontwisted) of $X$. For any $q\in X$, choose an orbifold chart $U_q/G_q$ satisfying
        Definition 3.1. A component of $f^{-1}(U_q/G_q)$ is of the form $U_q/\Gamma'$ for
    $\Gamma'\subset G_q$.
    It is clear that $G_q/\Gamma'$ is a subgroup of the orbifold fundamental
    group. Therefore, we obtain a group homomorphism
    $$\phi_q: G_q\rightarrow
\pi^{orb}_1(X).\leqno(4.1)$$
    It is easy to check that a different choice of component of $f^{-1}(U_q/G_q)$ or a
    different choice of $q\in X_{(g)}$ induces a
    homomorphism differing by a conjugation. Therefore, there is a unique
    map from the conjugacy classes of $G_q$ to the conjugacy classes of $\pi^{orb}_1(X)$.
    \vskip 0.1in
    \noindent
    {\bf Definition 4.1: }{\it We call $X_{(g)}$ a dormant sector if $\phi_p(g)=1$.}
    \vskip 0.1in
    If $X_{(g)}$ is a dormant sector, we define $L_{(g)}=1$. It will not
    receive any correction from discrete torsion. Non-dormant sectors are of
    the form $Y_g/C(g)$, where $Y_g\neq \emptyset$ is the fixed point loci of $1\neq g\in \pi^{orb}_1(X)$.
     $Y_{g}$ is a smooth suborbifold of $Y$.
    It is clear that $Y_{h^{-1}gh}$ is diffeomorphic to $Y_{g}$ by
    the action of $h$. By abusing the notation, we denote the twisted sector $Y_{g}/C(g)$ by
    $X_{(g)}$,
     where $C(g)$ is the centralizer of $g$.

    Let $\alpha$ be a global discrete
    torsion. Suppose that $S$ is the connected component of singular loci containing the image of $X_{(g)}$
     in $X$. We choose a small open
    neighborhood $U(S)$ of $S$ and suppose that local discrete torsion is $\alpha_S$.
    We replace $X$ by $U(S)$ in above construction.  We can use $L^{\alpha}_g$ to
    define a flat complex orbifold line-bundle
    $$L_{g}=Y_{g}\times_{L^{\alpha}_g}\C$$
    over $X_{(g)}$.
        \vskip 0.1in
    \noindent
    {\bf Lemma 4.2: }{\it
    \begin{description}
    \item[(1)] $L_{tgt^{-1}}$ is isomorphic to $L_{g}$ by the map
            $$t\times Id: Y_{g}\times \C\rightarrow Y_{tgt^{-1}}\times
                \C.\leqno(4.3)$$
                Hence, we can denote $L_{g}$ by $L_{(g)}$.
    \item[(2)] $L_{(g)}^{-1}=L_{(g^{-1})}.$
    \item[(3)] When we restrict to $X_{(g_1,\cdots,
    g_k)}=Y_{g_1}\cap\cdots \cap
                Y_{g_k}/C(g_1,\cdots, g_k)$,
            $L_{(g_1,\cdots, g_k)}=L_{(g_1)}\cdots L_{(g_k)}$, where
            $L_{(g_1,\cdots,  g_k)}=Y_{g_1}\cap\cdots
            Y_{g_k}\times_{\gamma_{g_1\cdots g_k}} \C.$
    \end{description}}
    \vskip 0.1in
    \noindent
    {\bf Proof: } Recall that there is an isomorphism
    $$t_{\#}: C(g)\rightarrow C(tgt^{-1})$$
    given by $t_{\#}(h)=tht^{-1}$. The map
    $$t:Y_{g}\rightarrow X_{tgt^{-1}}$$ is $t_{\#}$-equivariant. By Lemma 3.8,
    $\gamma_{tgt^{-1}}(tht^{-1})=\gamma_{g}(h)$ for $h\in C(g)$. Then,
    $$(t\times Id)(hx, \gamma(h)(v))=(thx, \gamma_{g}(h)(v))=(tht^{-1}tx,
    \gamma_{tgt^{-1}}(tht^{-1})(v)).\leqno(4.4)$$
    Then we take the quotient by $C(g), C(tgt^{-1})$ respectively to get
    an isomorphism between $L_{g}, L_{tgt^{-1}}$.
    (2) and (3) follow from the fact that for any $h\in C(g_1,\cdots,g_k)$,
    $$\gamma(\alpha)_{g_1\cdots g_k, h}=\gamma(\alpha)^{-1}_{h,g_1\cdots g_k}
    =\gamma(\alpha)^{-1}_{h,g_1}\cdots \gamma(\alpha)^{-1}_{h, g_k}=
        \gamma(\alpha)_{g_1, h}\cdots\gamma(\alpha)_{g_k, h}.\leqno(4.5)$$

    \vskip 0.1in
    \noindent
    {\bf Theorem 4.3: }{\it $\L_{\alpha}=\{L_{(g)}\}_{(g)\in T_1}$ is an inner local system of $X$.}
    \vskip 0.1in
    \noindent
    {\bf Proof: } Property (1) is obvious. The property (2) follows from Lemma 4.2.
    Let's prove property (3).
         Consider the
    image $\g'=(g'_1, g'_2, g'_3)$ of $\g$ in $\pi^{orb}_1(X)$ under the homomorphism (4.1). Then, we still have $g'_1g'_2g'_3=1$. There are three
    possibilities, (i) $g'_1=g'_2=g'_3=1$ and there is nothing to prove in
    this case; (ii) $g'_3=1, g'_2=(g')^{-1}_1$ is nontrivial; (iii)
    $g'_1,g'_2, g'_3$ are all nontrivial.

    For the second case, let $g=g'_1$. We have the following factorization
    $$e_1\times e_2\times e_3: X_{(\g)}\rightarrow X_{(g_1,g_2)}\times X_{(g_3)}
    \rightarrow X_{(g_1)}\times X_{(g_2)}\times X_{(g_3)}.$$
    However, $X_{(g_1, g_2)}=Y_{g}\cap Y_{(g^{-1})}/C(g,g^{-1})=Y_{g}/C(g).$
     Moreover, over $X_{g_1,g_2}$
    $$e^*_1 L_{(g)}e^*_2 L_{(g^{-1})}=L_{(g)}I^*L_{(g^{-1})}=1.\leqno(4.6)$$

    In the third case,
    $X_{(\g)}=Y_{g_1}\cap Y_{g_2}\cap Y_{g_3}/C(g_1, g_2, g_3)$. The
    proof follows from Lemma 4.2 (3).
    \vskip 0.1in
    \noindent
    {\bf Definition 4.4: }{\it Suppose that $\alpha$ is a global discrete torsion.
    We define the twisted orbifold cohomology $H^*_{orb, \alpha}(X, \C)=H^*_{orb}(X,
    \L_{\alpha}).$}
    \vskip 0.1in

    \section{Examples}

Only a few examples of global quotients have been computed by physicists \cite{VW}
    \cite{D}. It is still a very important problem to develop general machinery
    to compute discrete torsion and twisted orbifold cohomology. Here we compute
    five examples. First two have nontrivial discrete torsion. One is a global quotient and
    another one is a non-global quotient. The second example has the phenomenon that the most
    of twisted sectors are dormant sectors. The third one is Joyce example, where there is
    no nontrivial discrete torsion. However, there are nontrivial local systems. We will compute
    twisted orbifold cohomology given by  nontrivial local systems to match
    Joyce's desingularizations. Orbifold cohomology is strongly
    intertwine with group theory. We demonstrate it in last two
    examples.
    \vskip 0.1in
    \noindent
    {\bf Example 5.1 $T^4/\Z_2\times \Z_2$: } Here, $T^4=\C^2/\wedge$, where $\wedge$ is
    the lattice of integral points. Suppose that $g,h$ are generators of
    the first and the second factor of $\Z_2\times \Z_2$.
        The action of $\Z_2\times \Z_2$ on $T^4$ is defined as
    $$g(z_1, z_2)=(-z_1, z_2), h(z_1, z_2)=(z_1, -z_2).\leqno(5.1)$$
    The fixed point locus of $g$ is 4 copies of $T^2$. When we divide it by the remaining action
    generated by $h$, we obtain twisted sectors consisting of 4 copies of $S^2$. The
    degree shifting number for these twisted sectors is $\frac{1}{2}$. For the same reason,
    the fixed point locus of $h$ give twisted sectors consisting of $4$ copies of
    $S^2$ with degree shifting number $\frac{1}{2}$. The fixed point locus of $gh$ is
    16 points, which are fixed by the whole group. The degree shifting number of the 16 points
    is $1$. An easy calculation shows that nontwisted sector contributes one generator to degree
    0, 4 orbifold cohomology and two generators to degree 2 orbifold cohomology and no other.
    Using this information, we can compute the ordinary orbifold cohomology group
    $$b_{orb}^0=b_{orb}^4=1, b_{orb}^1=b_{orb}^3=8, b_{orb}^2=18.\leqno(5.2)$$

    By example 2.10, $H^2(\Z_2\times \Z_2, U(1))=\Z_2$. By Remark 2.2, the nontrivial generator of
    $H^2(\Z_2\times \Z_2, U(1))$ induces a discrete torsion $\alpha$. Next, we compute the
    twisted orbifold cohomology $H^*_{orb, \alpha}(T^4/\Z_2\times \Z_2, \C)$. Note that
    $\gamma (\alpha)_{gh,g}=\gamma (\alpha)_{gh,h}=-1$. Hence, the flat orbifold-bundles over
    the twisted sectors given by 16 fixed points of $gh$ are nontrivial. Therefore, they
    contribute nothing to twisted orbifold cohomology.
    For  two dimensional twisted sectors, let's consider a
    component of fixed point locus of $g$. By the previous description, it is $T^2$. $h$ acts on
    $T^2$. Then the twisted sector $S^2=T^2/\{h\}$. We observe that the flat orbifold line bundle
    over $S^2$ is constructed as $L=T^2\times_{L^{\alpha}_g} \C$. Hence $H^*(S^2, L)$
    is isomorphic to the space of invariant cohomology of $T^2$ under the action    of $h$ twisted by $\gamma(\alpha)_g$ as $h(\beta)=\gamma(\alpha)_{g,h}h^*\beta$.
    By example 2.10, $\gamma(\alpha)_{g,h}=-1$. The invariant cohomology is $H^1(T^2, \C)$. Using the degree
    shifting number to shift up its degree, we obtain the twisted orbifold cohomology
    $$b_{orb,\alpha}^0=b_{orb,\alpha}^4=1, b_{orb,\alpha}^1=b_{orb,\alpha}^3=0,     b_{orb,\alpha}^2=18.\leqno(5.3)$$
    \vskip 0.1in
    \noindent
    {\bf Example 5.2 $WP(2, 2d_1)\times WP(2, 2d_2)$ ($d_1,d_2>1, (d_1,d_2)=1$): } Here, $WP(2, 2d)$ is the weighted projective
    space of weighted $(2, 2d)$. $WP(2, 2d_1)\times WP(2, 2d_2)$ is not a global quotient unless
    $d_1=d_2=1$. In fact, its orbifold universal cover is $WP(1, d_1)\times WP(1, d_2)$ and
    $WP(2, 2d_1)\times WP(2, 2d_2)=WP(1, d_1)\times WP(1, d_2)/\Z_2\times \Z_2$. Hence, the orbifold
    fundamental group is $\Z_2\times \Z_2$. Therefore, there is
    a nontrivial discrete torsion $\alpha\in H^2(\Z_2\times \Z_2, U(1))$.

Next, we describe the
    twisted sectors. Suppose that $p=[0,1],q=[1,0]\in WP(1,d_1)$. We also
    use $p,q$ to denote its image in $WP(2,2d_1)$.  We use $p',q'$ to denote
    the corresponding points in $WP(1,d_2), WP(2,2d_2)$. $\{p\}\times WP(2,2d_2), \{p'\}\times
    WP(2,2d_1)$ give rise to two twisted sectors with degree shifting number $\frac{1}{2}$.
    $\{q\}\times WP(2,2d_2),\{q'\}\times WP(2,2d_1)$ give rise to $2d_1-1,2d_2-1$ many twisted sectors
    with degree shifting numbers $\frac{i}{2d_1},\frac{j}{2d_2}$ for $1\leq i\leq
    2d_1-1,1\leq j\leq 2d_2-1$.
    $\{p\}\times\{p'\}$ give rise to a twisted sector with degree shifting number $1$.
    $\{p\}\times \{q'\}$ give rise to $2d_2-1$-many twisted sectors with degree shifting numbers
    $\frac{1}{2}+\frac{i}{2d_2}$ for $1\leq i\leq 2d_2-1$. $\{q\}\times \{p'\}$
    give rise to $2d_1-1$-many twisted
    sectors with degree shifting numbers
    $\frac{1}{2}+\frac{i}{2d_1}$ for $1\leq i\leq 2d_1-1$. $\{q\}\times \{q'\}$
    give rise to $4d_1d_2-1$-many
    twisted sectors with degree shifting numbers $\frac{i}{2d_1}+\frac{j}{2d_2}$ for all $i,j$ except
    $(i,j)=(0,0)$. Using this information, we can write down ordinary orbifold cohomology
    $$b_{orb}^0=b_{orb}^4=1, b_{orb}^1=b^3_{orb}=6, b^2_{orb}=6$$
    $$b^{\frac{i}{d_1}}_{orb}=b^{\frac{i}{d_2}}_{orb}=1,
    b^{1+\frac{i}{d_1}}_{orb}=b^{1+\frac{i}{d_2}}_{orb}=3,
    b^{2+\frac{i}{d_1}}_{orb}=b^{2+\frac{i}{d_2}}_{orb}=2, \
    1\leq i\leq d_1-1, 1\leq j\leq d_2-1$$
    $$b^{\frac{i}{d_1}+\frac{j}{d_2}}_{orb}=1, 0\leq i\leq 2d_1-1, 0\leq j\leq 2d_2,
    (i,j)\neq (0,0),(d_1,d_2).
    \leqno(5.4)$$
    Next, we compute $H^*_{orb, \alpha}$. In this example, the most of twisted sectors are
    dormant sectors. To find nondormant sectors, recall that $WP(2,2d_1)\times WP(2,2d_2)=WP(1,d_1)\times WP(1,d_2)/Z_2\times Z_2$. Let $g$ be the
    generator of the first factor and $h$ be the generator of the second factor.
    The fixed points of $g$ is $\{p,q\}\times WP(1,d_2)$. We have two nondormant sectors
    obtained by modulo the
    remaining action generated by $h$. However, $\gamma(\alpha)_{g,h}=-1$.
    There is no invariant cohomology of $WP(1,d_2)$ under the
    action of $h$ twisted by $L^{\alpha}_g$. Hence, these two nondormant twisted sectors
    give no contribution to twisted orbifold cohomology. Their degree shifting numbers are 1.
    For the
    same reason, $WP(1,d_1)\times \{p',q'\}/g$ gives no
    contribution to twisted orbifold cohomology. The fixed point
    locus of $gh$ consists of  4 points which give 4 nondormant sectors. Again, their degree shifting
    numbers are 1. As we saw in last example, their flat orbifold bundles are nontrivial.
    Hence, they give no contribution to twisted orbifold cohomology.  Therefore,
    the twisted orbifold cohomology is
    $$b_{orb,\alpha}^0=b_{orb,\alpha}^4=1, b_{orb,\alpha}^1=b^3_{orb,\alpha}=2,
    b^2_{orb,\alpha}=2$$
    $$b^{\frac{i}{d_1}}_{orb,\alpha}=b^{\frac{i}{d_2}}_{orb,\alpha}=1,
    b^{1+\frac{i}{d_1}}_{orb,\alpha}=b^{1+\frac{i}{d_2}}_{orb,\alpha}=3,
    b^{2+\frac{i}{d_1}}_{orb}=b^{2+\frac{i}{d_2}}_{orb,\alpha}=2,
    \ 1\leq i\leq d_1-1, 1\leq j\leq d_2-1$$
    $$b^{\frac{i}{d_1}+\frac{j}{d_2}}_{orb,\alpha}=1,\ 0\leq i\leq 2d_1-1, 0\leq j\leq 2d_2,
    (i,j)\neq (0,0),(d_1,d_2).
    \leqno(5.5)$$
    \vskip 0.1in
    \noindent
    {\bf Example 5.3 $T^6/\Z_4$: } Here, $T^6=\C^3/\wedge$, where $\wedge$ is the lattice of
    integral points. The generator of $\Z_4$ acts on $T^6$ as
    $$\kappa: (z_1, z_2, z_3)\rightarrow (-z_1, iz_2, iz_3).\leqno(5.6)$$
    This example has been studied by D. Joyce \cite{JO}, where he constructed five different
    desingularizations. However, there is no  discrete torsion in the case which induces
    nontrivial orbifold cohomology.

    First all, the nontwisted sector contributes one generator to
    $H^{0,0}_{orb}, H^{3,3}_{orb}$, 5 generators to $H^{1,1}_{orb},
    H^{2,2}_{orb}$ and 2 generator to $H^{2,1}_{orb},
    H^{1,2}_{orb}$
    The fixed point loci of $\kappa, \kappa^3$ are 16-points
    $$\{(z_1, z_2, z_3)+\wedge: z_1\in \{0,\frac{1}{2}, \frac{i}{2}, \frac{1}{2}+\frac{i}{2}\},
    z_2, z_3\in \{0, \frac{1}{2}+\frac{i}{2}\}.$$
    These points are fixed by $Z_4$. Therefore, they generate
    32-twisted sectors in which 16 corresponds to the conjugacy
    class $(\kappa)$ and 16 corresponds to the conjugacy class
    $(\kappa^3)$. The sector with conjugacy class $(\kappa)$ has
    degree shifting number 1. The sector with conjugacy class
    $(\kappa^3)$ has degree shifting number 2.

    The fixed point loci of $\kappa^2$ is 16 copies of $T^2$,
    given by
    $$\{(z_1, z_2, z_3)+\wedge: z_1\in \C, z_2, z_3\in \{0,\frac{1}{2},
    \frac{i}{2}, \frac{1}{2}+\frac{i}{2}\}\}$$
    Twelve of the 16-copies of $T^2$ fixed by $\kappa^2$ are
    identified in pairs by the action of $\kappa$, and these
    contribute 6 copies of $T^2$ to the singular set of $T^6/Z_4$.
    On the remaining 4 copies $\kappa$ acts as $-1$, so these contribute
    4 copies of $S^2=T^2/\{\pm 1\}$ to singular set. The degree shifting number
     of these 2-dimensional twisted sectors  is 1.

     Next, we construct inner local systems. We start with two dimensional twisted sectors.
     Since $\kappa^{-2}=\kappa^2$, the condition (2) of Definition 2.1 tells us
    that the flat orbifold line bundle $L$ over two dimensional sectors has the
    property $L^2=1$. Now, we assign trivial line bundle to all
     $T^2$-sectors and $k (k=0,1,2,3,4)$-many $S^2=T^2/\{\pm 1\}$-sectors.
     For the remaining $S^2=T^2/\{\pm 1\}$-sectors, we assign a flat orbifold line
     bundle $T^2\times \C/\{\pm 1\}$. For the zero dimensional sectors, they are all
    points of two dimensional sectors. If we assign a trivial bundle on a two dimensional
    sector, we just assign trivial bundle to its point sectors. For these two
    dimensional sectors with nontrivial flat line bundle, we need to be careful to choose
    the flat orbifold line bundle on its point sectors to ensure the condition (3) of
    Definition 2.1. Suppose that $\Sigma$ is one of 2-dimensional sectors supporting
    nontrivial flat orbifold line bundle. It contains 4 singular points which generate the
    point sectors. Let $x$ be one of 4-points. $x$ generates two sectors given by the conjugacy
    classes $(\kappa), (\kappa^3)$. For condition (3),  we have to consider
     the conjugacy class of triple $(g_1,g_2, g_3)$ with
     $g_1g_2g_3=1$. The only nontrivial choices are
     $(\g)=(\kappa,\kappa, \kappa^2), (\kappa^2, \kappa^3,\kappa^3)$.
    The corresponding components of $X_{(\g)}$ are exactly the these singular points.
      Clearly, $x$ is a fixed by the whole group $Z_4$. The orbifold line bundle is given
    by the action of $\Z_4$ on $\C$. Consider the component of $X_{(\g)}$ generated by $x$.
    The pull-back of flat orbifold line bundle from 2-dimensional sector ($(\kappa^2)$-sector)
    is given by the action
    $\kappa v=-1$. A moment of thought tells us that for sectors $(\kappa), (\kappa^3)$, we should
    assign a flat
    orbifold line bundle given by the action of $\Z_4$ on $\C$ as $\kappa v=iv.$ It is easy to check
    that for above choices the condition (3) is satisfied for $X_{(\g)}$. Therefore,
    the twisted sectors given by $(x, (\kappa)), (x, (\kappa^3))$ give no contribution
    to twisted orbifold cohomology.  Suppose that the resulting local system is
     $\L_k$. For the sectors with trivial
     line bundle, they contribute $6+k$ generators to $H^{1,1}_{orb},
     H^{2,2}_{orb}$ and 6 generators to $H^{2,1}_{orb},
     H^{1,2}_{orb}$. Its point sectors contribute $4k$ generators to
    $H^{1,1}_{orb}, H^{2,2}_{orb}$.
     The remaining sectors contribute $4-k$ generators to
     $H^{2,1}_{orb}, H^{1,2}_{orb}$. Its point sectors give no contribution.
    Moreover, the nontwisted sector contributes
    $$h^{0,0}=h^{3,3}=2, h^{1,1}=5.$$
    In summary, we obtain
     $$\dim H^{0,0}_{orb}(T^6/Z_4, \L_k)=\dim H^{3,3}_{orb}(T^6/Z_4,
     \L_k)=1, \dim H^{1,1}_{orb}(T^6/Z_4, \L_k)=\dim H^{2,2}_{orb}(T^6/Z_4,
     \L_k)=11+5k,$$
    $$\dim H^{1,2}_{orb}(T^6/Z_4, \L_k)= \dim H^{2,1}_{orb}(T^6/Z_4,
     \L_k)=12-k\leqno(5.7)$$
     Our calculation matches the betti numbers of Joyce's desingularizations.

     The orbifold cohomology ring of following examples have been computed in
     \cite{CR}. Here, we compute their twisted version.
    \vskip 0.1in
    \noindent
    {\bf Example 5.4: } Let's consider the case that $X$ is  a point
    with a trivial group action of $G$. Suppose that $\alpha\in H^2(G, U(1))$ is a
    discrete torsion. We want to compute $H^*_{orb, \alpha}(X, \C)$. The twisted sector $X_{(g)}$ is a point
    with a group $C(g)$. It is obvious that $H^0(X_{(g)}, L^{\alpha}_g)=0$ unless
    $L^{\alpha}_g=1$. Recall that a conjugacy class $(g)$ is $\alpha$-regular iff $L^{\alpha}_g=1$.
    Hence, only $\alpha$-regular class will contribute.
    Therefore, torbifold cohomology is
    generated by $\alpha$-regular conjugacy classes of elements of $G$. All the
    degree shifting numbers are zero. By the same argument as
    nontwisted case, as a ring, $H^*_{orb, \alpha}(X, \C)$ is
    the  center of twisted group algebra $\C_{\alpha}[G]$.

    \vskip 0.1in
    \noindent
    {\bf Example 5.5: } Suppose that $G\subset SL(n, \C)$ is a
    finite subgroup.  Then,
    $\C^n/G$ is an orbifold. Suppose that $\alpha\in H^2(G, U(1))$
    is a discrete torsion. For any $g\in G$, the fixed point set
    $X_g$ is a vector subspace and  $X_{(g)}=X_g/C(g)$. By the
    definition, $L_{(g)}=X_g\times \gamma(\alpha)_g \C$.
    Therefore, $H^*(X_{(g)}, L_{(g)})$ is the  subspace of
    $H^*(X_g, \C)$ invariant under twisted action of $C(g)$
    $$h\circ \beta=\gamma(\alpha)_g(h) h^*\beta\leqno(5.8)$$
    for any $h\in C(g), \beta\in H^*(X_g, \C)$. However,
    $H^i(X_g, \C)=0$ for $i\geq 1$. Moreover, if
    $\gamma(\alpha)_g$ is nontrivial, $H^0(X_g, L_{(g)})=0$.
    Therefore,
    $H^{p, q}_{orb}=0$ for $p\neq q$ and
    $H^{p,p}_{orb}$ is
    a vector space generated by conjugacy class of
    $\alpha$-regular elements $g$ with $\iota_{(g)}=p$.
    Therefore, we have a natural decomposition
    $$H^*_{orb, \alpha}(X, \C)=Z[\C_{\alpha}[G])=\sum_p H_p,\leqno(5.9)$$
    where $H_p$ is generated by conjugacy classes of
    $\alpha$-regular elements $g$ with $\iota_{(g)}=p$.
    The ring structure is also easy to describe. Let $x_{(g)}$ be
    generator corresponding to zero cohomology class of twisted
    sector $X_{(g)}$ such that $g$ is $\alpha$-regular. The cup product is
    completely same as the nontwisted case except we replace conjugacy class
    by $\alpha$-conjugacy class. Let me sketch the calculation.
    As we showed in previous example,
    the multiplication of conjugacy classes can be described in terms of
    center of twisted group algebra $Z(\C_{\alpha}[G])$. But we have further
    restrictions in this case.  It is clear
    $$X_{(h_1,h_2,(h_1h_2)^{-1})}=X_{h_1}\cap X_{h_2}/C(h_1, h_2).$$
    To have nonzero invariant, we require that
    $$\iota_{(h_1h_2)}=\iota_{(h_1)}+\iota_{(h_2)}.\leqno(5.10)$$
     Then, we need to compute
   $$\int_{X_{h_1}\cap X_{h_2}/C(h_1,
   h_2)}e^*_3(vol_c(X_{h_1h_2}/C(h_1 h_2)))\wedge e(E),\leqno(5.11)$$
   where $vol_c(X_{h_1h_2}/C(h_1 h_2))$ is the compact supported top form with volume one.
   However,
   $$X_{h_1}\cap X_{h_2}/\subset X_{h_1h_2}$$
   is a submanifold. (5.11) is zero unless
   $$X_{h_1}\cap X_{h_2}=X_{h_1h_2}.\leqno(5.12)$$
   In this case, we call $(h_1,h_2)$ transverse. In this case,
   it is clear that obstruction bundle is trivial. Suppose that $d_{h_1, h_2}$ is the
   order of finite cover $X_{h_1 h_2}/C(h_1, h_2)\rightarrow X_{h_1 h_2}/C(h_1 h_2)$. Then,
   the  integral is $d_{h_1, h_2}$. Let
   $$I_{g_1, g_2}=\{(h_1, h_2); h_i \in (g_i),
   \iota_{(h_1)}+\iota_{(h_2)}=\iota_{(h_1h_2)}, (h_1,
   h_2)-transverse, (h_1h_2)-\alpha -regular \}.\leqno(5.13)$$
   Then,
   $$x_{(g_1)}\cup x_{(g_2)}=\sum_{(h_1,h_2)\in I_{g_1, g_2}}
   d_{h_1,h_2} x_{(h_1h_2)}.\leqno(5.14)$$

         \end{document}